\documentclass[12pt,a4paper]{article}

\usepackage[english]{babel}
\usepackage{amssymb}
\usepackage[latin1]{inputenc}
\usepackage{graphicx}
\usepackage{amsthm}
\usepackage{color}
\usepackage{stmaryrd}
\usepackage{fullpage}

\title{Transverse nonlinear instability of Euler-Korteweg solitons}
\date{February 2015}
\author{\textsc{Matthew Paddick}\footnote{Sorbonne Universit\'es, UPMC Univ Paris 06, UMR 7598, Laboratoire Jacques-Louis Lions, F-75005, Paris, France}~\footnote{CNRS, UMR 7598, Laboratoire Jacques-Louis Lions, F-75005, Paris, France} \\ \texttt{paddick@ljll.math.upmc.fr}}

\begin{document}

\renewcommand{\labelitemi}{$\bullet$}
\newtheorem{theo}{Theorem}[section]
\newtheorem{lemma}[theo]{Lemma}
\newtheorem*{nota}{Notation}
\newtheorem*{assu}{Assumption}
\newtheorem{propo}[theo]{Proposition}
\newtheorem{coro}[theo]{Corollary}
\newcommand{\preu}[1]{\textit{\underline{#1}}}
\newcommand{\thref}[1]{Theorem \ref{#1}}
\newcommand{\lemref}[1]{Lemma \ref{#1}}
\newcommand{\propref}[1]{Proposition \ref{#1}}
\newcommand{\cororef}[1]{Corollary \ref{#1}}
\newcommand{\tworef}[3]{#3 \ref{#1} and \ref{#2}}

\newcommand{\til}[1]{\tilde{#1}}
\newcommand{\rplus}{\mathbb{R}^{+}}
\newcommand{\rtwo}{\mathbb{R}^{2}}
\renewcommand{\div}{\mathrm{div}~}
\newcommand{\divt}{\mathrm{div}\tang}
\newcommand{\rot}{\overrightarrow{\mathrm{rot}~}}
\newcommand{\eps}{\varepsilon}
\newcommand{\norme}[2]{\left\| #2 \right\| _{#1}}
\newcommand{\tilnorme}[2]{|| #2 ||_{#1}}
\newcommand{\trinorme}[2]{\interleave #2 \interleave _{#1}}
\newcommand{\tang}{_\tau}
\newcommand{\derp}[1]{\partial_{#1}}
\newcommand{\diag}{\mathrm{diag}}
\newcommand{\dOmega}{\partial\Omega}
\newcommand{\Lip}{\mathrm{Lip}}
\newcommand{\rhs}[1]{\mathbf{RHS}(\ref{#1})}
\newcommand{\mathand}{~~\mathrm{and}~~}
\newcommand{\mathif}{\mathrm{if}~}
\newcommand{\sdel}[1]{#1_{\delta}}
\newcommand{\vphi}{\varphi}
\newcommand{\bb}[1]{\mathbb{#1}}
\newcommand{\real}{\mathrm{Re}}
\newcommand{\mattwo}[4]{\left( \begin{array}{cc} #1 & #2 \\ #3 & #4 \end{array} \right)}
\newcommand{\id}{\mathrm{Id}}
\newcommand{\ap}{^{\mathrm{ap}}}

\maketitle

\begin{abstract}
We show that solitary waves for the 2D Euler-Korteweg model for capillary fluids display nonlinear instability when subjected to transverse perturbations.
\newline

\newline

\end{abstract}

\section{Introduction}

We consider the motion of a compressible, inviscid and isentropic planar fluid, in which internal capillarity is taken into account. This phenomenon occurs for example at diffuse interfaces in liquid-vapour mixes \cite{BDDJ}.
 In this model, the free energy of the fluid depends on both the density of the fluid, the scalar function $\rho$, and its gradient $\nabla \rho$ in the following way:
$$ F(\rho,\nabla\rho) = F_0(\rho) + \frac{1}{2} K(\rho) |\nabla \rho|^2, $$
with $K$ and $F_0$ two given smooth, positive functions. We then derive the pressure from the free energy as follows:
$$ P(\rho,\nabla \rho) = \rho \frac{\derp{}F}{\derp{}\rho} - F = P_0(\rho)  + \frac{1}{2}(\rho K'(\rho)-K(\rho))|\nabla \rho|^2, $$
in which $P_0$ is the standard part of the pressure. Let $g_0(\rho)$ be the bulk chemical potential of the fluid, so that $\rho g'_0(\rho)=P_0'(\rho)$.
 Then, the principles of classical mechanics yield the Euler-Korteweg equation that we will study:
\begin{equation} \left\{ \begin{array}{rcl} \derp{t}\rho + \div(\rho u) & = & 0 \\ \derp{t}u + (u \cdot\nabla) u & = & \nabla \left(K(\rho)\Delta \rho + \frac{1}{2} K'(\rho) |\nabla \rho|^2 - g_0 (\rho)\right) \\ (\rho,u)|_{t=0} & = & (\rho_0,u_0). \end{array} \right. \label{ek} \end{equation}
The variables are $t\in\bb{R}^+$ and $(x,y)\in\bb{R}^2$; as is standard, the operators $\nabla$, $\div$and $\Delta$ contain only derivatives with respect to the space variables $x$ and $y$.
 The unknowns of equation (\ref{ek}) are the density $\rho$ and the velocity vector field $u:\bb{R}^+\times \bb{R}^2 \rightarrow \bb{R}^2$. The scalar functions $g_0$ and $K$ are given.
\newline

In this paper, we will be interested in the transverse stability of solitary wave solutions of (\ref{ek}). These are 1D travelling waves written as
$$ Q_c(t,x) = \left( \begin{array}{c} \rho_c(t,x) \\ u_{c}(t,x) \end{array} \right) = Q_c(x-ct), $$
with $u_c$ scalar (not a 2D vector field).
 Based on a remark by T. Benjamin \cite{Bt}, S. Benzoni, R. Danchin, S. Descombes and D. Jamet showed in \cite{BDDJ} that the hamiltonian structure of the system
 led to the existence of solitary wave solutions for every $c\in\bb{R}$ in the case of a van der Waals-type pressure law.
 There are two kinds of travelling-wave solutions, depending on the homoclinic or heteroclinic nature of the connecting orbits.
 In the homoclinic case, the wave has identical endstates, and we can write
$$ \lim_{|z|\rightarrow +\infty} Q_c(z) = Q_{c,\infty} = (\rho_{c,\infty}, u_{c,\infty}), $$
and such a travelling wave solution is called a soliton. We will be studying this type of solitary wave.
 In the heteroclinic case, the endstates are different, but must nonetheless satisfy a Rankine-Hugoniot-type condition. See \cite{BDDJ} for more on these solutions, which are called kinks.

From now on, we set $c$ and $Q_c$ is a soliton, whose endstate, which we also set (thus we drop the index $c$), satisfies
\begin{equation} \rho_\infty g'_0(\rho_\infty) > (u_\infty-c)^2 , \label{saddle} \end{equation}
which means that $Q_\infty$ is a saddle point for the hamiltonian ODE satisfied by $Q_c$. Under this condition, we have that $\rho_c'$ vanishes only once.
\newline

The standard Lyapunov stability notion is that if solutions of equation (\ref{ek}) have initial conditions close to $Q_c$, they remain close to $Q_c(x-ct)$ at all times.
 But this notion is not satisfactory in describing the stability of travelling waves. Indeed, let $c'\neq c$ be close to $c$; we then have that the profile $Q_{c'}$ is close to $Q_c$, but, as the speeds are different, $Q_{c'}(x-c't)$ and $Q_c(x-ct)$ drift apart, despite their profiles remaining very similar.
 To see this, for a given $t$, compare $Q_{c'}(x-c't)$ with the translated profile $Q_c (x-c't)=Q_c(x-ct+(c-c')t)$.
 The correct notion of stability therefore stems from taking the difference of solutions with all the translated versions of $Q_c(x-ct)$.
 A travelling wave solution will be considered stable if it is \underline{orbitally stable}: for every $\eps>0$, there exists $\delta>0$ such that if $\norme{}{(\rho_0,u_0)-Q_c}\leq \delta$, then
$$ \sup_{t\in\bb{R}^+} \inf_{a\in\bb{R}} \norme{}{(\rho(t),u(t)) - Q_c(\cdot-a-ct)} \leq \eps. $$

The problem of orbital stability can be divided into two parts, depending on the type of perturbation we consider.
 1D perturbations are perturbations of $Q_c$ that depend only on $x$ and satisfy $u_2(t,x)=0$. The stability problem associated with these perturbations has been in part dealt with by S. Benzoni \textit{et al.} in \cite{BDDJ}, and improved upon by Benzoni in \cite{Bs}.
 A sufficient condition for orbital stability was obtained in the first paper using an argument by M. Grillakis, J. Shatah and W. Strauss \cite{GSS}, while the second article adds a sufficient condition for linear instability.
 See also J. H\"owing \cite{Hj2, Hj} for other stability results for the 1D Euler-Korteweg system.

The question of transverse stability deals with perturbations that also depend on the transverse variable $y$ and have a 2D velocity field.
 So far, Benzoni in \cite{Bs} and F. Rousset and N. Tzvetkov in \cite{RT-MRL} have proved linear instability. This occurs when the linearised equation around $Q_c$ has eigenvalues with positive real part.
 On one hand, Benzoni used Evans functions computations to get that orbitally stable solitons are transversally linearly unstable.
 On the other hand, Rousset and Tzvetkov applied a simple, abstract criterion for instability in linearised PDEs with a hamiltonian structure in the case where the endstate of the soliton satisfies (\ref{saddle}).
 We recall this criterion in \thref{mrl}. This criterion was applied to other equations with solitary waves in the same article: KP-I and Gross-Pitaevskii.
\newline

The result of this paper is that the spectral instability mentioned above implies nonlinear instability of Euler-Korteweg solitons.
\begin{theo} Let $Q_c(x-ct)$ be a soliton solution to (\ref{ek}) such that the endstate $Q_\infty$ satisfies (\ref{saddle}).
 Then there exist $\delta_0$ and $\eps_0>0$ such that for every $0<\eps<\eps_0$, there exists an initial condition $U_0=(\rho_0,u_0)$ with
$$ \norme{H^{s}(\bb{R}^2)}{U_{0}-Q_c}\leq \eps $$
for some $s>0$, such that, for every $a\in\bb{R}$, the solution $U=(\rho,u)$ of (\ref{ek}) with this initial condition satisfies, at a time $T^\eps\sim \ln(\eps^{-1})$,
$$ \norme{L^2(\bb{R}^2)}{U(T^\eps)-Q_c(\cdot-cT^\eps-a)} \geq \delta_0 . $$ 
Moreover, the velocity $u$ can be chosen to be potential: $u=\nabla \varphi$. \label{ek-inst} \end{theo}

The proof relies on an argument originally by E. Grenier \cite{Ge}, in which one constructs an approximate solution $U\ap$ to the equation based on a WKB expansion starting with the reference solution and unstable eigenmodes.
 An energy estimate on $U-U\ap$ shows that, if there are enough terms in the expansion, this difference is small, thus, for times under $T^\eps$, the linear instability is dominant.
 Primarily used to obtain nonlinear instability of boundary layers in numerous settings (unstable Euler shear flows and Prandtl layers \cite{Ge}, Ekman layers for rotating fluids \cite{DG, Rf, MRekman}, Ekman-Hartman layers in MHD \cite{DG}, Navier-Stokes with a boundary-layer-scale slip condition \cite{Pm11}),
 the method relies on building approximate solutions from a compact set of unstable wave numbers, thus the idea has been transposed to showing transverse nonlinear instability of solitary waves, when these can be shown to be linearly unstable.
 F. Rousset and N. Tzvetkov have thus obtained nonlinear instability of solitary waves in many models:
 KP-I and NLS \cite{RT1}, multiple hamiltonian models including generalised KP-I and the Boussinesq equations \cite{RT3}, and the free-surface water-waves equation \cite{RT2}.
\newline

Similarly to many situations, such as KP-I, our result extends to the $y$-periodic framework $\bb{R}\times\bb{T}_L$, where $\bb{T}_L$ is the torus of length $L>0$, in the following way:
 there exists a critical period $L_0(c)>0$, depending on the soliton's speed, such that if $L>L_0(c)$, we have \thref{ek-inst}.
 The proof is identical, and limited to above a critical period due to the loss of linear instability below this period (the set of unstable wavenumbers is bounded).
\newline

\textbf{Outline of the proof.} The proof of \thref{ek-inst} is in two parts.
 First, in section 2, we build on Rousset and Tzvetkov's linear instability theorem, deriving an important resolvent estimate for the Euler-Korteweg equation, linearised around $Q_c$.
 This will allow us, in section 3, to build an approximate solution $U\ap$ with the appropriate behaviour of being predominantly unstable for $t \sim T^\eps$.
 An energy estimate on $U-U\ap$ will then be used. This is particularly important ensure that the time of existence of $U$ is large enough to get the desired amplification,
 as we only have local existence for solutions to (\ref{ek}) (see S. Benzoni, R. Danchin and S. Descombes in \cite{BDD}).
 Article \cite{BDD} also provides a blow-up criterion, but the instability phenomenon is not related to the blow-up if it occurs;
 indeed, the mechanism is also observed on systems that have global solutions (most of the previous examples, whether they concern boundary layer or solitary wave instability, fall in this category).
 Finally, combining the two will lead to the instability result.

\section{Linear analysis}

In this part, we consider the linearised problem about $(\rho_c,u_c)$ to obtain an essential ingredient in order to apply Grenier's method.

Considering that $u$ is potential, we write the system on $(\rho,\vphi)$, where $u=\nabla \vphi$:
$$ \left\{ \begin{array}{rcl} \derp{t}\rho + \nabla \vphi\cdot\nabla \rho+ \rho\Delta\vphi & = & 0 \\ \derp{t}\vphi + \frac{1}{2}|\nabla \vphi|^2 & = & K(\rho)\Delta\rho + \frac{1}{2} K'(\rho)|\nabla\rho|^2 - g_0(\rho) , \end{array} \right. $$
which we linearise around $(\rho_c,u_c)$. Having changed the space variable from $x$ to $x-ct$ (which turns the solitary wave into a stationary solution), we are interested in
\begin{equation} \left\{ \begin{array}{rcl} \derp{t}\rho & = & (c\derp{x}-u_c\derp{x} -u_c')\rho - (\rho_c' \derp{x} + \rho_c\Delta) \vphi \\ \derp{t}\vphi & = & (c\derp{x}-u_c\derp{x})\vphi + (K(\rho_c)\Delta + K'(\rho_c)\rho_c' \derp{x} - m)\rho , \end{array} \right. \label{ekpl} \end{equation}
with $m = K'(\rho_c)\rho_c''+K'(\rho_c)(\rho_c')^2-g_0'(\rho_c)$. We abbreviate the system by defining two operators
$$ J = \mattwo{0}{1}{-1}{0} \mathand {\cal L} = \mattwo{-\derp{x}(K(\rho_c)\derp{x}) - K(\rho_c)\derp{yy}^2 - m}{(u_c-c)\derp{x}}{-\derp{x}((u_c-c)\cdot)}{-\derp{x}(\rho_c\derp{x})-\rho_c\derp{yy}^2} , $$
thus the system (\ref{ekpl}) can be written synthetically as $\derp{t}V=J{\cal L}V$, where $V=(\rho,\vphi)$.
\newline

The first part of the linear analysis involves finding unstable eigenmodes for (\ref{ekpl}).
 These are non-trivial solutions to the equation that can be written as $V(t,x,y)=e^{\sigma t}e^{iky}v(x)$ for $k\neq 0$ and $\real(\sigma)>0$.
 We rewrite (\ref{ekpl}) using the Fourier transform on the transverse variable $y$: the equation $\derp{t}V=J{\cal L}V$ becomes $\sigma v = JL(k)v$ with
$$ L(k) = \mattwo{-\derp{x}(K(\rho_c)\derp{x}) + K(\rho_c)k^2 - m}{(u_c-c)\derp{x}}{-\derp{x}((u_c-c)\cdot)}{-\derp{x}(\rho_c\derp{x})+\rho_c k^2} . $$

We begin by examining the existence of eigenmodes and the behaviour of $\sigma$ depending on $k$, and we follow up with an important resolvent estimate for $JL(k)$.

\begin{propo} \underline{Properties of the linearised equation.}
\newline

(a) The linearised equation is unstable, that is there exist eigenmodes written as $V(t,x,y)=e^{\sigma t}e^{iky}v(x)$, with $v\in H^2(\bb{R})$ and $\real(\sigma)>0$, that solve (\ref{ekpl}).
 For each $k$, the dimension of the subspace of unstable solutions of $\sigma v=JL(k)v$ is at most 1.

The instability is localised in the transverse Fourier space: there exists $k_{\max}>0$ such that, for $|k|\geq k_{\max}$, eigenvalues necessarily satisfy $\real(\sigma)\leq 0$.
 Let $\sigma(k)$ be the eigenvalue of $JL(k)$ with highest real part. Then the function $k\mapsto \real(\sigma(k))$ has a global maximum $\sigma_0>0$ at a certain $k_0>0$.
\newline

(b) If $V(t,x,y) = e^{iky}U(t,x)$, we define the following semi-norm for $U$:
$$ \norme{X^j_k}{U(t)}^2 = \norme{H^{j+1}(\bb{R})}{U_1(t)}^2 + \norme{H^j(\bb{R})}{\derp{x}U_2(t)}^2 + |k|^2\norme{H^j(\bb{R})}{U(t)}^2 . $$
It is essentially the $H^{j}$ norm of $|k|U(t)$ plus the $H^{j+1}$ norm of $U(t)$, omitting the $L^2$ norm of $U_2(t)$.

Set $\gamma>\sigma_0$, $n>0$, $s\in\bb{N}$ and let $U$ solve
\begin{equation} \derp{t}U(t,x)=JL(k)U(t,x)+F(t,x,k), \label{linpb}\end{equation}
with $U(0,x)=0$, and, $F$ satisfying, for every $j\leq s$,
\begin{equation} \norme{H^{j+1}(\bb{R})}{\derp{t}^{s-j}F(k,t)} \leq M_s \frac{e^{\gamma t}}{(1+t)^n} , \label{fbound} \end{equation}
uniformly for $|k|\leq k_{\max}$ ($M_s$ does not depend on $k$). Then $U$ satisfies a similar estimate: for every $j\leq s$, and uniformly for $|k|\leq k_{\max}$, we have
\begin{equation} \norme{X^j_k}{\derp{t}^{s-j}U(t)} \leq C_s (1+t)^{-n} e^{\gamma t} . \label{resolvent} \end{equation}
\label{resest} \end{propo}

\textbf{\underline{Consequences of part (b):}} a quick energy estimate on the equation of $U_2$ yields that $U_2(t)\in L^2(\bb{R}^2)$ (as $\vphi(0)=0$), and this $L^2$ norm also satisfies (\ref{resolvent}).
 We will therefore subsequently consider that the result is valid in $H^s$, for any $s\geq 0$.

By the Parseval equality, this result also implies identical $H^s$ bounds for $(\rho,u)$ when written as $(\rho,u)=(\rho,\nabla\vphi)$ and $(\rho, \vphi) = \int_{\bb{R}} f(y)e^{iky} U(t,x)~dk$, with $f\in{\cal C}^{\infty}_0(\bb{R})$.
 Indeed, norms of $|k|^2 U$ can be replaced, using equation (\ref{linpb}), by derivatives on $x$ and $t$ that satisfy (\ref{resolvent}).

\subsection{Proof of \propref{resest} (a) \\ Properties of eigenmodes}

\subsubsection{Existence of unstable eigenmodes.}

The existence of unstable eigenmodes was shown by F. Rousset and N. Tzvetkov \cite{RT-MRL} using a general criterion for detecting transverse linear instability of solitary waves.
 We have seen that equation (\ref{ekpl}) for functions written as $V(t,x,y)=e^{\sigma t}e^{iky} v(x)$ becomes an eigenvalue problem, that is
\begin{equation} \sigma v=JL(k)v \label{eigpb} \end{equation}
with $J$ a skew-symmetric matrix and ${\cal L}$ a self-adjoint differential operator on $(L^2(\bb{R}^2))^2$ whose domain is seen to be $(H^2(\bb{R}^2))^2$. We have the following result for such systems.

\begin{theo} (Rousset and Tzvetkov, \cite{RT-MRL})

If $L$ has the following properties,
\begin{itemize}
\item (H1) there exists $k_{\max}>0$ and $\alpha>0$ such that $L(k)\geq \alpha \id$ for $|k|\geq k_{\max}$;
\item (H2) for every $k\neq 0$, the essential spectrum of $L(k)$ is included in $[\alpha_k,+\infty[$ with $\alpha_k>0$;
\item (H3) $L'(k)$ is a positive operator;
\item (H4) the spectrum of $L(0)$ consists of one isolated negative eigenvalue $-\lambda$ and a subset of $\rplus$;
\end{itemize}
then there exist $\sigma>0$ and $k\neq 0$ such that (\ref{eigpb}) has a non-trivial solution, and, for every unstable wavenumber $k$, such an eigenvalue $\sigma$ is unique.
\label{mrl} \end{theo}

This is shown by finding $k'>0$ such that $L(k')$ has a one-dimensional kernel, and by using the Lyapunov-Schmidt method in the vicinity of this point; we do not detail the proof of this theorem.
 Proof that the linearised Euler-Korteweg system satisfies the hypotheses of this theorem was also done in \cite{RT-MRL}, but we shall briefly recall this, as it contains some useful arguments for the subsequent points of \propref{resest} (a).
\newline

(H1): using Young's inequality, $ab\leq \frac{\delta}{2}a^2 + \frac{1}{2\delta} b^2$, with $\delta=\frac{K(\rho_c)}{2}$, we quickly get that
\begin{equation} (L(k)v,v) \geq \int_{\bb{R}} \frac{K(\rho_c)}{2} |\derp{x}v_1|^2 + \left(k^2-\frac{1}{2}\right) |v_1|^2 + \rho_c |\derp{x}v_2|^2 + \left(k^2-\frac{1}{K(\rho_c)}\right)|v_2|^2, \label{H1} \end{equation}
which is greater than $\alpha\norme{L^2}{U}^2$ for $|k|$ large enough (remember that $\rho_c$ and $K(\rho_c)$ are positive).
\newline

(H2): as $\lim_{x\rightarrow \pm\infty} Q_c(x) = Q_\infty=(\rho_\infty,u_\infty)$, with standard arguments \cite{Hd}, and using the fact that $L(k)$ is self-adjoint, the essential spectrum of $L(k)$ is given by that of
$$ L_\infty(k) = \mattwo{K(\rho_\infty)(-\derp{xx}^2+k^2) + g'_0(\rho_\infty)}{(u_\infty-c)\derp{x}}{(c-u_\infty)\derp{x}}{\rho_\infty(-\derp{xx}^2+k^2)} , $$
whose essential spectrum can be determined by using the Fourier transform in the $x$-variable and explicitly writing the eigenvalues $\mu(\xi,k)$.
 We get that these are positive when $k\neq 0$. The essential spectrum of $L(k)$ is equal to that of $L_\infty(k)$, so (H2) is verified.
\newline

(H3): we easily have $L'(k)=\mathrm{diag}(2kK(\rho_c),2k\rho_c)$.
\newline

(H4): we apply the following lemma to $L(0)$.
\begin{lemma} Let $L$ be a symmetric operator on a Hilbert space such that
$$ L = \mattwo{L_1}{A}{A^*}{L_2} $$
with $L_2$ invertible. Then, we can write
$$ (Lv,v) = \left((L-AL_2^{-1}A^*)v_1,v_1\right) + \left( L_2(v_2+L_2^{-1}A^* v_1), v_2+L_2^{-1} A^* v_1\right) . $$ \label{alglem} \end{lemma}
As a result, we write
$$ (L(0)v,v) = (Mv_1, v_1) + \int_{\bb{R}} \rho_c \left| \derp{x}v_2 + \frac{1}{\rho_c}(u_c-c) v_1\right|^2 ~dx, $$
with $M=-\derp{x}(K(\rho_c)\derp{x}\cdot)-m-\frac{(u_c-c)^2}{\rho_c}$. $M$ is a second-order differential operator on which we can perform Sturm-Liouville analysis \cite{DS}.
 First, the essential spectrum of $M$ is included in $[\alpha,+\infty[$ with $\alpha>0$; indeed $M$ is a perturbation of $M_\infty = -K(\rho_\infty)\derp{xx}^2+g_0 '(\rho_\infty)-\frac{(u_\infty-c)^2}{\rho_\infty}$, whose essential spectrum is positive under the assumption that $\rho_\infty g'_0(\rho_\infty) > (u_\infty-c)^2$.
 Next, the function $\rho_c'$ is in the kernel of $M$, and it has one zero, so by Sturm-Liouville theory, $M$ has a unique negative eigenvalue associated with an eigenfunction $U_1^-$. Setting $U_2^-$ such that $\derp{x}v_2^- = \frac{-1}{\rho_c}(u_c-c)v_1^-$, we have a generalised eigenfunction for $L(0)$ (the second component is not in $L^2$).
 By using $H^2$ approximations of $U_2^-$, we see that $(L(0)v,v)$ can be negative with $U\in H^2$, confirming that $L(0)$ has one negative eigenvalue.
\newline

\subsubsection{Localisation of instability and boundedness of unstable eigenvalues.}

We start be taking the real part of the $L^2$ scalar product of the eigenvalue equation (\ref{eigpb}) by $L(k)U$: we get that $ \real(\sigma) (L(k)U,U) = 0$, as $J$ is skew-symmetric.
 The operator $L(k)$ satisfies (H1) of \thref{mrl}, so, if $\real(\sigma)>0$, we must have $0=(L(k)U,U)\geq \alpha\norme{L^2}{U}^2$ for $|k|\geq k_{\max}$.
 Thus, the only function satisfying $\sigma U = JL(k)U$ with $\real(\sigma)>0$ and $|k|$ large is $U=0$; there are no unstable eigenfunctions for $|k|$ large.
\newline

In order to get the boundedness of the unstable eigenvalues, we decompose $L(k)$ as follows: $L(k) = L_0(k) + L_1$ with
$$ L_0(k):=\mattwo{-\derp{x}(K(\rho_c)\derp{x})+K(\rho_c)k^2-m_0}{0}{0}{-\derp{x}(\rho_c \derp{x}) + \rho_c k^2} , $$
where $-m_0(x) = \max\left(-m(x),\frac{1}{2}g_0'(\rho_\infty)\right)$. We compute the scalar product of (\ref{eigpb}) and $L_0(k)v$, and take the real part, which gives us
\begin{equation} \real(\sigma)(L_0(k)v,v) = \real(JL_1 v, L_0(k)v). \label{jl1} \end{equation}
It is quickly noticed that there exists $\alpha>0$ such that
\begin{equation} \real(\sigma)(L_0(k)v,v) \geq \alpha\real(\sigma)\left( \norme{L^2}{\derp{x}v}^2 + k^2\norme{L^2}{v}^2 + \norme{L^2}{v_1}^2 \right). \label{l0low} \end{equation}
We shall now bound $|\real(JL_1 v, L_0(k) v)|$ by the same norms as on the right. Note that
$$ JL_1 = \mattwo{-\derp{x}((u_c-c)\cdot)}{0}{m-m_0}{-(u_c-c)\derp{x}} . $$
Computing the scalar product directly, and bounding the terms involving $\rho_c$ and $u_c$ in $L^\infty$, we get
\begin{eqnarray*} |\real(JL_1 v,L_0(k)v)| & \leq & C\left[(1+k^2)(|\real(v_1,\derp{x}v_1)|+|\real(\derp{x}v_2,v_2)|) + |\real(v_1,\derp{x}v_2)| \right. \\
 & & \hspace{20pt} + \norme{L^2}{\derp{x}v_1}^2 + k^2\norme{L^2}{v_1}^2 + \norme{L^2}{\derp{x}v_2}^2 \\
 & & \hspace{20pt} \left. +|(v_1,\derp{xx}v_1)| + |(v_1,\derp{xx}^2 v_2)| + |(|k|v_1,|k| v_2)| \right] \end{eqnarray*}
The top line vanishes by integration by parts, and, also by integrating by parts, the first two terms in the final line are repeats of terms in the second line.
 Finally, by using Young's inequality, the last term is bounded by $k^2 \norme{L^2}{v}^2$, thus proving that, combining with (\ref{jl1}) and (\ref{l0low}), there exists $C>0$ such that
$$ \real(\sigma)\left(\norme{L^2}{\derp{x}v}^2 + k^2\norme{L^2}{v}^2 + \norme{L^2}{v_1}^2 \right) \leq C \left(\norme{L^2}{\derp{x}v}^2 + k^2\norme{L^2}{v}^2 + \norme{L^2}{v_1}^2 \right), $$
which implies that $\real(\sigma)$ cannot be unbounded when positive. This ends the proof of part (a) of \propref{resest}.

\subsection{Proof of \propref{resest} (b) \\ Resolvent estimate}

The proof of part (b) is split in two. First, we get the result for $s=0$; we bound the $X^0_k$ norm of $U$ by similar norms of $F$ by using the Laplace transform and spectral arguments.
 The case $s>0$ is then obtained by induction on $s$, the number of total derivatives (time and space).

\subsubsection{The case $s=0$}

The proof of (\ref{resolvent}) for $s=0$ relies on the Laplace transform, and is similar to the resolvent estimate proofs in \cite{RT1} and \cite{Pm11}.
 Let $\sigma_0<\gamma_0<\gamma$. For $f(t)$, we denote by $\tilde{f}(\tau)$ the following Laplace transform:
$$ \tilde{f}(\tau) := \int_0^{+\infty} \exp(-(\gamma_0+i\tau)t) f(t)~dt. $$
Using the Laplace transform turns equation (\ref{linpb}), $\derp{t}U = JL(k)U + F$, into an eigenvalue problem:
\begin{equation} (\gamma_0 + i\tau) \tilde{U}(\tau) = JL(k) \tilde{U}(\tau) + \tilde{F}(\tau). \label{lappb} \end{equation}
As $\gamma_0>\sigma_0$, $\gamma_0+i\tau$ is not in the spectrum of $JL(k)$.
 Indeed, we can use the strategy employed to prove that hypothesis (H2) is satisfied to show that the essential spectrum of $JL(k)$ is embedded in $i\bb{R}$.
 Once again using the argument from \cite{Hd}, we can examine the spectrum of the Fourier transform in $x$ of $JL_\infty(k)$,
$$ {\cal F}_x(JL_\infty)(\xi,k) = \mattwo{-i(u_\infty-c)\xi}{\rho_\infty(\xi^2+k^2)}{-K(\rho_\infty)(\xi^2+k^2)-g_0'(\rho_\infty)}{-i(u_\infty-c)\xi}, $$
which contains the solutions of the equation
$$ X^2 - 2i\xi(u_\infty-c)X + \rho_\infty K(\rho_\infty)(\xi^2+k^2)^2 + \rho_\infty g'_0(\rho_\infty)k^2 + (\rho_\infty g'_0(\rho_\infty)-(u_\infty-c)^2)\xi^2 = 0, $$
which depend on $(\xi,k)$. Using the positiveness of $\rho_\infty$, $K(\rho_\infty)$ and condition (\ref{saddle}), we get that the discriminant of this equation is negative for $(\xi,k)\neq 0$,
 and clearly the only eigenvalue at $(\xi,k)=(0,0)$ is zero, so the essential spectrum of $JL(k)$ is imaginary.
\newline

As $\gamma_0+i\tau$ is not in the spectrum of $JL(k)$ for any $\tau\in\bb{R}$, the norm of the resolvent $((\gamma_0+i\tau)\id-JL(k))^{-1}$ is uniformly bounded for $(\tau,k)$ in any compact subset of $\bb{R}^2$.
 It remains to show that, for $|k|\leq k_{\max}$, there exists the following bound for $|\tau|$ large.
\begin{lemma} If $\tilde{U}$ solves (\ref{lappb}), then there exists $C,~M>0$ such that, for $|\tau|\geq M$,
\begin{equation} \tilnorme{X^0_k}{\tilde{U}(\tau)} \leq C \tilnorme{H^1}{\tilde{F}(\tau)}. \label{x0est} \end{equation} \label{estlem} \end{lemma}

\preu{Proof:} we consider the scalar product of the Laplace-transformed equation (\ref{lappb}) with $L(k)\tilde{U}$, and write
\begin{equation} (\gamma_0+i\tau)(L(k)\tilde{U},\tilde{U}) = (\tilde{F},L(k)\tilde{U}). \label{splap} \end{equation}
Note that $(L(k)\tilde{U},\tilde{U}) = (L(0)\tilde{U},\tilde{U}) + K(\rho_c)k^2\tilnorme{L^2}{\tilde{U}_1}^2+\rho_c k^2\tilnorme{L^2}{\tilde{U}_2}^2$. Let us concentrate on the term $(L(0)\tilde{U},\tilde{U})$. Using \lemref{alglem}, we know that it is equal to
$$ (L(0)\tilde{U},\tilde{U}) = (M\tilde{U}_1,\tilde{U}_1) + \int_{\bb{R}} \rho_c \left|\derp{x}\tilde{U}_2+\frac{1}{\rho_c}(u_c-c)\tilde{U}_1\right|^2~dx, $$
The operator $M$, which we remind the reader is equal to $-\derp{x}(K(\rho_c)\derp{x}\cdot)-m-\frac{(u_c-c)^2}{\rho_c}$, and whose quadratic form is defined on $H^1$, has one simple negative eigenvalue, as well as a one-dimensional kernel containing $\rho_c'$.

Recall $U^-$ the generalised eigenfunction corresponding to the negative eigenvalue of $L(0)$. We do not have $U_2^-\in L^2$, so we consider $U^-=(U^-_1,0)$. We denote $U^0=(\rho_c',0)$, which is in the kernel of $M\otimes \id$. We renormalise $U^-$ and $U^0$ so that their $L^2$ norms are equal to 1. Let $U^+$ be orthogonal to $U^-$ and $U^0$; we show that
\begin{equation} (L(0)U^+,U^+)\geq \eta \norme{X^0_0}{U^+} \label{perpest} \end{equation}
for some $\eta>0$. First, as $U^+_1$ is not in the kernel or the negative eigenspace of $M$, we have $(MU^+_1,U^+_1)\geq \alpha \norme{H^1}{U^+_1}^2$, since the essential spectrum of $M$ is included in $[\alpha,+\infty[$. Thus, we already have
\begin{equation} (L(0)U^+,U^+) \geq \alpha \norme{H^1}{U^+_1}^2. \label{mpos} \end{equation}
But this does not suffice to get the $X^0_0$ norm. Using the fact that $\rho_c$ is positive, there is a positive $\beta$ such that
\begin{eqnarray*} \int_{\bb{R}} \rho_c \left| \derp{x}U^+_2+\frac{1}{\rho_c}(u_c-c)U^+_1\right|^2 ~dx & \geq & \beta \norme{L^2}{\derp{x}U^+_2 + \frac{1}{\rho_c}(u_c-c)U^+_1}^2 \\
 & \geq & \beta \norme{L^2}{\derp{x}U^+_2}^2 + C\beta \norme{L^2}{U^+_1}^2 - 2C|(\derp{x}U^+_2,U^+_1)| \end{eqnarray*}
for some $C>0$. We use Young's inequality on the final term with $\delta=\beta/2C$, thus there exists $C'\in\bb{R}$ such that
$$ (L(0)U^+,U^+) \geq (\alpha+C')\norme{H^1}{U^+_1}^2 + \frac{\beta}{2} \norme{L^2}{\derp{x}U^+_2}^2. $$
If perchance $C'$ is negative, we add $\frac{|C'|}{\alpha} \times(\ref{mpos})$ to the above, and obtain that there does indeed exist $\eta>0$ such that we have (\ref{perpest}), and $(L(k)U^+,U^+) \geq \eta \norme{X^0_k}{U^+}^2$.
\newline

We now write the orthogonal decomposition in $L^2$ of the first component, $\tilde{U}_1 = aU^-_1 + bU^0_1 + U^+_1$, and replace in (\ref{splap}).
 The eigenfunctions $U^-$ and $U^0$ are fixed, so their $H^1$ norms are given constants. On the right-hand side, using integration by parts and basic estimates including Young's inequality, we have
$$ |(\tilde{F},L(k)\tilde{U})| \leq  \frac{\eta}{4} \norme{X^0_k}{U^+}^2 + C(\tilnorme{H^1}{\tilde{F}}^2+a^2+b^2), $$
while on the left-hand side, we have
$$ (L(k)\tilde{U},\tilde{U}) \geq \eta (|k|^2\tilnorme{L^2}{U^+}^2 + \norme{X^0_0}{U^+}^2) - C(a^2 + (|a|+|b|)\norme{X^0_0}{U^+}^2). $$
Taking the real part of (\ref{splap}) and moving the negative part of the above to the right-hand side and once again applying Young's inequality to absorb $\norme{X^0_0}{U^+}$, we get
\begin{equation} \norme{X^0_k}{U^+}^2 \leq C(\tilnorme{H^1}{\tilde{F}}^2 + a^2 + b^2 ). \label{endperp} \end{equation}

To finish off, we take the dot product of (\ref{lappb}) with $U^-$ and $U^0$. We quickly get
$$ (\gamma_0+i\tau) a = -(\tilde{U},L(k)JU^-) + (\tilde{F},U^-) \mathand (\gamma_0+i\tau)b = -(\tilde{U},L(k)JU^0) + (\tilde{F},U^0) . $$
Integrating by parts as usual, we get
$$ a^2+b^2 \leq \frac{C}{\gamma_0^2 + |\tau|^2}(a^2+b^2+\norme{X^0_k}{U^+}^2 + \tilnorme{H^1}{\tilde{F}}^2) .$$
We see that if $|\tau|$ is large enough, $a^2+b^2$ can be absorbed by the left-hand side. Combining with (\ref{endperp}), we get the result. $\square$
\newline

For the end of the proof of \propref{resest} (b), we start by using the Parseval equality in the following,
$$ \int_0^T e^{-2\gamma_0 t} \norme{X^0_k}{U(t)}^2~dt = \int_0^T \tilnorme{X^0_k}{\tilde{U}(t)}^2~dt \leq C\int_0^T \tilnorme{H^1}{\tilde{F}(t)}~dt = C\int_0^T e^{-2\gamma_0 t}\norme{H^1}{F(t)}~dt. $$
Now recall the assumption on $F$, (\ref{fbound}): we have
$$ \int_0^T e^{-2\gamma_0 t} \norme{X^0_k}{U(t)}^2~dt \leq CM_0 \int_0^T \frac{e^{2(\gamma-\gamma_0)t}}{(1+t)^\alpha}~dt \leq C_0 \frac{e^{2(\gamma-\gamma_0)T}}{(1+T)^n}. $$
We inject this in the energy estimate on (\ref{linpb}), that is
$$ \frac{d}{dt} (\norme{X^0_k}{U(t)}^2) \leq C(\norme{X^0_k}{U(t)}^2 + \norme{H^1}{F(t)}^2), $$
and multiply the result by $e^{-2\gamma_0 t}$, integrate in time and we get the result. $\square$

\subsubsection{The induction for $s>0$}

The extension of \propref{resest} (b) to every $s\geq 0$ will follow the lines of the similar result on the linearised water-waves equation in \cite{RT2};
 it is done with a double induction, double in the sense that one is embedded in the other.

The first induction is on $s$, the total number of derivatives. Set $s>0$, and we assume that, for every $s'<s$ and $j\leq s'$, we have (\ref{resolvent}), that is
$$ \norme{X^j_k}{\derp{t}^{s'-j} U(t)} \leq C_{s'} \frac{e^{\gamma t}}{(1+t)^n}. $$
To get the wanted result, we must prove that, for every $0\leq j\leq s$,
\begin{equation} \norme{\dot{X}_k}{U}^2 := \norme{\dot{H}^1}{\derp{t}^{s-j}\derp{x}^j U}^2 + |k|^2 \norme{L^2}{\derp{t}^{s-j}\derp{x}^j U}^2 \leq C_s \frac{e^{2\gamma t}}{(1+t)^{2n}}, \label{dotx} \end{equation}
where $\dot{H}^1$ is the usual homogeneous Sobolev norm on $\bb{R}$. This is done by induction on the number of space derivatives, $j$.
 Note that the $\dot{X}_k$ norm (semi-norm if $k=0$) defined here is a sort of homogeneous Sobolev norm expressed in the Fourier space.
 Morally, at rank $s$ of the induction, we must get bounds for the $L^2$ norms of terms involving $s+1$ derivatives.
\newline

Starting with $j=0$, we are interested in the $X^0_k$ norm of $\derp{t}^{s}U$.
 Simply differentiate the equation, (\ref{linpb}), $s$ times with respect to $t$, and notice that $W(t)=\derp{t}^s U(t)-\derp{t}^{s}U(0)$ satisfies $\derp{t}W = JL(k)W + G$, with $W|_{t=0}=0$ and $\norme{H^{s+1}}{G(t)}\leq 2M_s (1+t)^{-n}e^{\gamma t}$: we can re-use the case $s=0$ shown above.
\newline

Now, let $j>0$. To lighten the notations, we will write $V_{s,j}=\derp{t}^{s-j}\derp{x}^j V$.
 Note that we want to control the $\dot{X}_k$ norm of $U_{s,j}$, which means $s+1$ derivatives in total, of which $j+1$ space or Fourier derivatives.
 We apply $\derp{t}^{s-j}\derp{x}^j$ to the equation. This time, the derivatives do not commute with $JL(k)$, hence we consider
$$ \derp{t} U_{s,j} = JL(k) U_{s,j} + J [\derp{x}^j,L(k)] U_{s-j,0} + F_{s,j} := JM_{s,j}(k)U + F_{s,j}. $$
We take the real part of the scalar product of this equation with $M_{s,j}(k)U$, which yields
$$ \frac{1}{2} \frac{d}{dt} (U_{s,j},L(k)U_{s,j}) = -\real(U_{s+1,j},[\derp{x}^j,L(k)]U_{s-j,0}) + \real (F_{s,j},M_{s,j}(k)U). $$
To bound the second part of the right-hand side, we look more closely at the commutator term in $M(k)U$. We notice that there exist two sets of $L^\infty$ matrices $(m^1_i,m^2_i)_{0\leq i\leq j+1}$ such that
\begin{equation} [\derp{x}^j,L(k)] \derp{t}^{s-j} U = \left(\sum_{i=0}^{j+1} m^1_i(x) U_{s-j+i,i} \right) + \left(\sum_{i=0}^{j-1} m_i^2(x) k^2 U_{s-j+i,i} \right). \label{commu} \end{equation}
We notice that all the terms, \textit{except} the one with $i=j+1$, have a total of $s$ derivatives or less, and thus, using $k^2\leq k_{\max}|k|$, they are controlled by our induction hypothesis.
 Integrating the terms of $(F_{s,j},L(k)U_{s,j})$ involving $j+2$ space derivatives and using assumption (\ref{fbound}) and Young's inequality with a parameter $\eta$ to be chosen later, we obtain that the right-hand side is bounded by
$$ |(F_{s,j},M_{s,j}(k)U)| \leq \frac{\eta}{2}\norme{\dot{X}_k}{U_{s,j}}^2 + C\frac{e^{2\gamma t}}{(1+t)^{2n}} . $$
To deal with the first term of the right-hand side, we notice that $U_{s+1,j}$ has $s+1$ derivatives, of which $j$ space derivatives, hence the $L^2$ norm of $U_{s+1,j}$ falls under our induction hypothesis.
 We can thus use Young's inequality and use decomposition (\ref{commu}) once again, and get
$$ \frac{1}{2}\frac{d}{dt} (U_{s,j},L(k)U_{s,j}) \leq \eta\norme{\dot{X}_k}{U_{s,j}}^2 + C\frac{e^{2\gamma t}}{(1+t)^{2n}}. $$

Finally, we integrate this in time, and recall (\ref{H1}) from the verification of the (H1) hypothesis of \thref{mrl}, which says that
$$ (L(k)U_{s,j},U_{s,j}) \geq \eta' \norme{\dot{X}_k}{U_{s,j}}^2 - C\norme{L^2}{U_{s-1,j-1}}^2, $$
in which the final term can be moved to the right-hand side and controlled by the induction hypothesis.
 Note that we can take the $L^2$ norm of $U_{s-1,j-1}$ since the number of time derivatives is preserved (we do not encounter the second component of $U$ which is not assumed to be in $L^2$).
 In total, we therefore have
$$ \norme{\dot{X}_k}{U_{s,j}(T)}^2 \leq \frac{\eta}{\eta'} \int_0^T \norme{\dot{X}_k}{U_{s,j}(t)}^2~dt + C\frac{e^{2\gamma t}}{(1+t)^{2n}}. $$
We choose $\eta$ in the Young inequalities above so that $\eta/\eta' \leq 2\gamma$, and the Gr\"onwall lemma gives us (\ref{dotx}) for the couple $(s,j)$. Both inductions are now complete. $\square$

\section{Nonlinear instability}

In this part, $U=(\rho,u)$. Obtaining \thref{ek-inst} relies on the construction of an approximate solution $U\ap$ built around unstable eigenmodes of the linearised equation.
 In our case, this construction is classical and we will not write all the details of the calculations (see also, for instance, \cite{Ge, DG, RT3, RT2}).
 Energy estimates must then be obtained on $U-U\ap$ to ensure that the approximate solution is close enough to the exact solution for long enough to see the difference between $U\ap$ and $Q_c$ reach an amplitude ${\cal O}(1)$.
 This is the more delicate part, as we remind the reader that only local existence is guaranteed for solutions of the Euler-Korteweg equation, thus we must also ensure that the solution $U$ still exists when the instability appears.
 The closing argument is standard for this method.

\subsection{Construction and properties of the approximate solution}

For a whole number $N$ independent of $\eps$ to be chosen later, we will set
$$ U\ap = Q_c(t,x) + \sum_{j=1}^N \eps^j U_j(t,x,y) , $$
with $u_j$ potential: we set $V\ap = (\rho\ap, \vphi\ap)$, $V_j=(\rho_j,\vphi_j)$. This is expected to solve the Euler-Korteweg system leaving an error of order $\eps^{N+1}$, as follows,
$$ \left\{ \begin{array}{rcl} \derp{t}\rho\ap + \div (\rho\ap \nabla \vphi\ap) & = & \eps^{N+1}R\ap_1 \\
\derp{t}\vphi\ap + \frac{1}{2} |\nabla\vphi\ap|^2 + g_0(\rho\ap) & = & K(\rho\ap)\Delta\rho\ap + \frac{1}{2}K'(\rho\ap)|\nabla \rho\ap|^2 + \eps^{N+1}R\ap_2, \end{array} \right. $$
thus, expanding $V\ap$ in these equations and isolating the terms of order $\eps$, we see that $V_j$ solves the linearised Euler-Korteweg equation around $Q_c$ with a source term:
\begin{equation} \derp{t} V_j = J{\cal L}V_j + R_j, \label{levelj} \end{equation}
with $R_j$ containing nonlinear interaction terms between the $V_n$ with $n<j$, but with the sum on indices in each interaction term equal to $j$.
 For instance, we have $R_1=0$, and then, in the equation on $\rho_2$, we see that $R_{2,1}=\div (\rho_1 \nabla \vphi_1)$.
 In the equation on $\vphi_j$, we notice that there will be quadratic terms similar to the one we have mentioned, but also cubic interaction terms, in the sense that they involve $V_{n_1}$, $V_{n_2}$ and $V_{n_3}$ with $n_1+n_2+n_3=j$.
 Likewise, $R\ap$ contains all the interaction terms whose sum of indices is greater than $N$. We do not detail the remainders any further.
\newline

We now construct $V_1$ as a wavepacket of unstable eigenmodes of the linearised equation. Recall that $k_0>0$ is the global maximum for the function
$$ \tilde{\sigma}: k\mapsto \max \{\real(\lambda)~|~\lambda \in \sigma(JL(k))\} , $$
where $\sigma(JL(k))$ is the spectrum of the operator $JL(k)$. We assume that this maximum is nondegenerate (although the method is easily adapted to the degenerate case, see \cite{RT2}).
 The function $\tilde{\sigma}$ is continuous, hence we can define
$$ V_1(t,x,y) = \int_{\bb{R}} f_1(k) e^{iky} e^{\tilde{\sigma}(k)t} v_1(k,x)~dk , $$
with $f_1(k)$ smooth, even, equal to 1 in the vicinity of $k_0$ and supported in the set $\{k~|~\tilde{\sigma}(k)>3\sigma_0/4 \}$, and $w(t,x)=e^{\tilde{\sigma}(k)t} v_1(k,x)$ solving $\derp{t}w=JL(k)w$.
 For any $s\geq 0$, we consider the quantity
$$ \sum_{s'=0}^s \int_{\bb{R}} f_1^2(k)|k|^{2s'} C_{s,s'} e^{2\tilde{\sigma}(k)t} ~dk = \sum_{s'=0}^s \int_{\bb{R}} f_1^2(k) |k|^{2s'} \norme{H^{s-s'}(\bb{R})}{w(k,t)}^2~dk = \norme{H^s(\bb{R}^2)}{V_1(t)}^2 , $$
the second equality being the Parseval equality.
 We use the Laplace method around the critical point $k_0$ to get that $\displaystyle \norme{H^s}{V_1(t)}^2 \mathop{\sim}_{t\rightarrow +\infty} t^{-1/2} e^{2\sigma_0 t}$, so
\begin{equation} \frac{1}{C_{1,s}} \frac{e^{\sigma_0 t}}{(1+t)^{1/4}} \leq \norme{H^s(\bb{R}^2)}{V_1(t)} \leq C_{1,s} \frac{e^{\sigma_0 t}}{(1+t)^{1/4}}. \label{laplm} \end{equation}

We get estimates on $V_j$ by induction. Assume that $I=\mathrm{supp}(f_1)$ is made up of two separate intervals around $\pm k_0$, and we set
$$ V_j(t,x,y) = \int_I \cdots \int_I v_j(k_1,\cdots,k_j;t,x) e^{ik_1 y}\cdots e^{ik_j y}~dk_1\cdots dk_j. $$
Assuming that, for every $n<j$,
\begin{equation} \norme{H^s}{v_n(k_1,\cdots,k_i;t)} \leq C_n \exp[n(\tilde{\sigma}(k_1)+\cdots+\tilde{\sigma}(k_n))t] , \label{induc} \end{equation}
we get that $v_j$ solves the linearised Fourier-transformed equation
$$ \derp{t}v_j(k_1,\cdots,k_j) = JL(k_1+\cdots+k_j) v(k_1,\cdots,k_j) + r_j(k_1,\cdots,k_j), $$
in which $\norme{H^s}{r_j(k_1,\cdots,k_j)} \leq C\exp[n(\tilde{\sigma}(k_1)+\cdots+\tilde{\sigma}(k_j))t]$ by the structure of the remainder (combination of products between $v_n$'s with $n<j$) and (\ref{induc}).
 Then, since, for $k\in I$, $\tilde{\sigma}(k)>3\sigma_0/4$, the sum in the exponential is greater than $\sigma_0$, and setting $v_j|_{t=0}=0$, we can apply \propref{resest} (b) to get that $v_j(k_1,\cdots,k_j;t)$ satisfies (\ref{induc}).
 We then use Parseval's equality and the Taylor expansion of $\tilde{\sigma}$ around the critical point $k_0$ to write that, for some $\beta>0$,
\begin{eqnarray*} \norme{H^s}{V_j(t,x,y)}^2 & = & \int_{jk\in jI} \norme{H^s}{\int_{k_1+\cdots+k_j=jk} v_j(k_1,\cdots,k_j;t,x) e^{ijky}~dk_1\cdots dk_{j-1}}^2~dk \\
 & \leq & \int_{jk\in jI} Ce^{2(j\sigma_{0}-j\beta(k-k_{0})^{2})t}\int_{\sum_{m=1}^j k_{m}=jk} e^{-2\beta\sum_{m=1}^j (k_{m}-k)^{2}t}~dk_1\cdots dk_{j-1}~dk \end{eqnarray*}
Integrate these gaussian functions (remembering that $k_j=jk-\sum_{n=1}^{j-1} k_n$), and we get the desired inequality: for every $j\leq 1$,
\begin{equation} \norme{H^s(\bb{R}^2)}{V_j(t)} \leq C_j \frac{e^{j\sigma_0 t}}{(1+t)^{j/4}}. \label{hsuj} \end{equation}

We now take a look at the remainder of the equation on $V\ap$, $R\ap$, which contains the interaction terms of the equation whose sum of indices is greater than $N$. Similarly to our proof of (\ref{hsuj}), we have
\begin{equation} \norme{H^s}{\eps^{N+1}R\ap(t)} \leq \sum_{j=N+1}^{3N} C_j \frac{\eps^{j} e^{j\sigma_0 t}}{(1+t)^{j/4}}, \label{rap1} \end{equation}
and, in what follows, we will be interested in times for which the smaller powers of $\eps(1+t)^{-1/2}e^{\sigma_0 t}$ are dominant. We set $T^*_\eps$ such that
$$ \frac{\eps e^{\sigma_0 T^*_\eps}}{(1+T^*_\eps)^{1/4}} = \kappa, $$
for $0<\kappa<1$ to be chosen later, and $t=T^*_\eps-\tau$. Replace in (\ref{rap1}), and we have
\begin{eqnarray*} \norme{H^s}{\eps^{N+1} R\ap(T^*_\eps-\tau)} & \leq & \left(\max_{j\in\{N+1,\cdots,3N\}} C_j\right) \sum_{j=N+1}^{3N} \kappa^j e^{-j\sigma_0 \tau} \\
 & \leq & C_R \kappa^{N+1} e^{-(N+1)\sigma_0 \tau}, \end{eqnarray*}
which, returning to the original time variable $t$, gives us, for $t\leq T^*_\eps$,
\begin{equation} \norme{H^s}{\eps^{N+1}R\ap(t)} \leq C_R \frac{\eps^{N+1}e^{(N+1)\sigma_0 t}}{(1+t)^{(N+1)/4}}. \label{hsrap} \end{equation}

\subsection{Getting the instability}

If $U$ is the solution of the Euler-Korteweg system (\ref{ek}) with the initial condition $U(0)=U\ap(0)$, we will observe the instability by studying
$$ \norme{L^2(\bb{R}^2)}{U(t) - Q_c(t)} \geq \norme{L^2(\bb{R}^2)}{U\ap(t)-Q_c(t)} - \norme{L^2(\bb{R}^2)}{U(t)-U\ap(t)}. $$
On one hand, we have $U\ap-Q_c = \sum_{j=1}^N \eps^j U_j$, and
\begin{eqnarray*} \norme{L^2(\bb{R}^2)}{\sum_{j=1}^N \eps^j U_j} & \geq & \norme{L^2(\bb{R}^2)}{\eps U_1(t)} - \sum_{j=2}^N \norme{L^2(\bb{R}^2)}{\eps^j U_j} \\
 & \geq & C'_1 \frac{\eps e^{\sigma_0 t}}{(1+t)^{1/4}} - \sum_{j=2}^N C_j \frac{\eps^j e^{j\sigma_0 t}}{(1+t)^{j/4}} \end{eqnarray*}
by (\ref{laplm}) and (\ref{hsuj}). Taking times smaller than $T^*_\eps$, we can consider that the sum on the right behaves like $\eps^2(1+t)^{-1}e^{2\sigma_0 t}$, and, replacing $t$ by $T^*_\eps-\tau$, we write
\begin{eqnarray*} \norme{L^2}{(U\ap-Q_c)(T^*_\eps-\tau)} & \geq & \kappa \left[ C'_1 e^{-\sigma_0 \tau} - \kappa C'_2 e^{-2\sigma_0 \tau} \right] \\
 & \geq & \kappa C_1' e^{-\sigma_0\tau} ( 1 - \frac{\kappa C'_2}{C'_1} e^{-\sigma_0\tau} ). \end{eqnarray*}
We notice that, for a given $C>0$, there exists $\tau_C>0$ such that, for $\tau\geq \tau_C$, $1-Ce^{-\sigma_0\tau} \geq 1/2$, so we set $\tau_1>0$, independent of $\eps$, such that, for $\tau\geq \tau_1$,
\begin{equation} \norme{L^2(\bb{R}^2)}{(U\ap-Q_c)(T^*_\eps-\tau)} \geq \frac{\kappa C'_1}{2} e^{-\sigma_0\tau} . \label{uapmino} \end{equation}

On the other hand, we will use energy estimates to ensure that $\norme{L^2(\bb{R}^2)}{U(t)-U\ap(t)}$ is small.
 We will readily use those shown by S. Benzoni, R. Danchin and S. Descombes in \cite{BDD}, which are obtained by considering the equation on $(G,z)=(G, u + iw)$, with $G$ a primitive of the function $\rho\mapsto \sqrt{K(\rho)/\rho}$ and $w=\nabla (G(\rho))$.
 The equation on $z$ is a Schr\"odinger-type equation, written as
$$ \derp{t}z + u\cdot\nabla z + i\nabla z\cdot w + i\nabla ({\cal A}(\rho) \div z) = {\cal Q}(\rho) , $$
with functions ${\cal A}$ and ${\cal Q}$ that we do not detail, while $G$ satisfies
$$ \derp{t}G + (u\cdot\nabla) G + {\cal A}(\rho) \div u = 0 . $$
The approximate solution satisfies a similar system with a remainder term which we will denote ${\cal R}=({\cal R}_G,{\cal R}_z)$.
 From now on, we use tildes to designate the difference between the exact and approximate terms in this system, \textit{e.g.} $\tilde{u} = u-u\ap$, $\tilde{\cal A} = {\cal A}(\rho)-{\cal A}(\rho\ap)$.
 The difference $(\tilde{G},\tilde{z})$ satisfies the equation
\begin{equation} \left\{ \begin{array}{l} \derp{t}\tilde{G} + u\cdot\nabla \tilde{G} + \tilde{u}\cdot\nabla G\ap + {\cal A}(\rho)\div\tilde{u} + \tilde{A}\div u\ap = -{\cal R}_G \\
\derp{t}\tilde{z} + u\cdot\nabla \tilde{z} + \tilde{u}\cdot\nabla z\ap + i\nabla \tilde{z}\cdot w + i\nabla z\ap \cdot \tilde{w} \\
\hspace{65pt} + i\nabla({\cal A}(\rho)\div \tilde{z}) + i\nabla (\tilde{\cal A}\div z\ap) = \tilde{Q} - {\cal R}_z. \end{array} \right. \label{schro} \end{equation}
The energy estimates in $H^s$ on this equation involve multiplying by an adequately chosen gauge $\psi_s$.
 In the potential case, this gauge is ${\cal A}(\rho)^{s/2}$ (whereas in the non-potential case, a additional term is needed), and the weighted norm $\norme{L^2}{\psi_s \Lambda^s \cdot}$, in which $\Lambda^s$ is the standard Fourier multiplier for $s$ derivatives, is equivalent to the standard $H^s$ norm.
 We send the reader to article \cite{BDD} for details, in particular part 3 where the reasons for the simpler norm is explained, and part 6 where energy estimates on the difference between a solution and a reference solution ((\ref{schro}) without the source term ${\cal R}$) are obtained.
\newline

Using the $H^s$ norms of $Q_c-Q_\infty$ as constants, the energy estimate reads
$$ \frac{d}{dt}\norme{H^s}{\tilde{z}}^2 \leq C\norme{H^s}{\tilde{z}}(\tilnorme{L^2}{\tilde{G}}+\norme{H^s}{\tilde{z}})(1+\norme{H^{s+1}}{\nabla z\ap})(1+\norme{H^{s-1}}{\nabla z\ap}+\norme{H^s}{\tilde{z}}) + C\norme{H^s}{\cal R}^2. $$
We rewrite this as follows: let $W=U-U\ap$; there exists a polynomial function $Q$ such that, for $s>0$ large enough (to handle the $L^\infty$ norms by Sobolev embedding),
$$ \norme{H^s}{W(t)}^2 \leq \int_0^t Q(\norme{H^s}{U\ap(z)-Q_\infty}+\norme{H^s}{U(z)})\norme{H^s}{W(z)}^2 + \norme{H^s}{R\ap(z)}^2~dz. $$
We will now choose $N$ to get the right growth in time for $W$, as well as $\kappa$ to get the existence up to $T^*_\eps$ of the exact solution $U$.
 First of all, in the same way that we get (\ref{uapmino}), we note that
\begin{eqnarray*} \norme{H^s}{U\ap(t)-Q_\infty} & \leq & \norme{H^s}{Q-Q_\infty} + \sum_{j=1}^N C_j \frac{\eps^j e^{j\sigma_0 t}}{(1+t)^{j/4}} \\
 & \leq & \norme{H^s}{Q-Q_\infty} + \kappa \end{eqnarray*}
when $t\leq T^*_\eps-\tau_2$, with $\tau_2\geq \tau_1$ independent of $\eps$.
 We consider times $t\leq T_W$ so that $\norme{H^s}{W(t)}\leq 1$ and $\rho(t,x,y)>0$ (no vacuum on the exact solution), and choose $N$ so that, for $t\leq T_W$,
$$ 2N\sigma_0 > Q(\norme{H^s}{Q-Q_\infty}+1+\kappa) . $$
A variant of the Gr\"onwall inequality from \cite{PSV} then provides us with
\begin{equation} \norme{H^s}{W(t)}^2 \leq C \frac{\eps^{N+1} e^{2(N+1)\sigma_0 t}}{(1+t)^{(N+1)/2}} \label{uexmajo1} \end{equation}
for $t\leq T_W$. Now, take $t=T^*_\eps-\tau$: we notice that the right-hand side is smaller than $C(N)\kappa^{2(N+1)}$, which is therefore smaller than $\kappa$ if $\kappa<1$ is small enough.
 We now choose $\kappa$ so that $2\kappa < \min \rho_c$, to ensure that there is no vacuum. We therefore have $T_W\geq T^*_\eps$ by a bootstrap argument.
 So, (\ref{uexmajo1}) is valid for $t=T^*_\eps-\tau$ with $\tau\geq \tau_2$, and we have
\begin{equation} \norme{H^s}{W(T^*_\eps-\tau)} \leq C'_0 \kappa^{N+1} e^{-(N+1)\sigma_0 \tau}. \label{uexmajo} \end{equation}

Recall that $f_1$, which is the localising function for the unstable term $U_1$, is supported in $I$ which is made up of two closed intervals that do not contain 0.
 Set $f$ a smooth function such that $f(k)=1$ for $k\in I$ and $f=0$ in the vicinity of 0, and let us define $\Pi$, a Fourier projector on frequencies in $I$, by
$$ {\cal F}_y(\Pi u)(x,k) = f(k)({\cal F}_y u)(x,k). $$
As $U(0,x,y) = Q_c(x)+\eps U_1(0,x,y)$, $\Pi U|_{t=0} = U_0$. Moreover, for any $a\in\bb{R}$, the difference $Q_c(x-ct-a)-Q_c(x-ct)$ does not depend on $y$, hence $\Pi(Q_c(\cdot-ct-a)-Q_c(\cdot-ct)) = 0$.
 We can now combine (\ref{uapmino}) and (\ref{uexmajo}) to show the instability: we have, for any $a$,
\begin{eqnarray*} \norme{L^2(\bb{R}^2)}{U(t)-Q_c(\cdot-a-ct)}|_{t=T^*_\eps-\tau} & \geq & \norme{L^2(\bb{R}^2)}{\Pi(U(t)-Q_c(\cdot-ct)}|_{t=T^*_\eps-\tau} \\
& \geq & \frac{\kappa C'_1}{2}e^{-\sigma_0\tau} \left[ 1 - \frac{2\kappa^N C'_0}{C'_1} e^{-N\sigma_0 \tau} \right] , \end{eqnarray*}
For $\tau\geq \tau_3\geq \tau_2$, we have the last exponential on the right smaller than $1/2$, and, as a result, letting $\tau'\geq \tau_3$ be fixed, independent of $\eps$,
$$ \norme{L^2}{(U-Q_c)(T^*_\eps-\tau')} \geq \frac{\kappa C'_1}{4}e^{-\sigma_0 \tau'} := \delta_0 . $$
The number $\delta_0$ we have found is positive and does not depend on $\eps$: \thref{ek-inst} is proved. $\square$
\newline

\textbf{\underline{Acknowledgements:}} this result was obtained during my Ph.D at Rennes 1 University, under the supervision of Fr\'ed\'eric Rousset, who I warmly thank for the opportunity to work on the topic,
 and for his guidance throughout the preparation of this article.
 
The author is partially supported by the ANR project DYFICOLTI, Agence Nationale de la Recherche grant ANR-13-BS01-0003-01.

\begin{small}
\bibliography{transbib}
\bibliographystyle{abbrv}
\end{small}

\end{document}